\documentclass[12pt, leqno]{article}
\usepackage{amssymb}
\usepackage{amsmath}
\def\q \m#1#2{{\raise1pt\hbox{$#1$}\kern-1pt\big/
               \kern-1pt\raise-1pt\hbox{$#2$}}}

\def\bN{{\rm \bf N}}
\def\bR{{\rm \bf R}}
\def\bZ{{\rm \bf Z}}
\def\bQ{{\rm \bf Q}}
\def\bC{{\rm {\bf C}}}
\def\ch{{\rm  ch}}
\def\bH{{\rm \bf H}}
\def\sR{{ \rm \scriptsize  \bf R}}
\def\sZ{{ \rm \scriptsize  \bf Z}}
\def\sC{{ \rm \scriptsize  \bf C}}

\def\sQ{{ \rm \scriptsize  \bf Q}}

\def\Sym{{\rm  Sym}}

\newfam\msbfam
\font\twelmsb=msbm10 at 12pt
\font\tenmsb=msbm10 at 10 pt
\font\sevenmsb=msbm10 at 7pt

\textfont\msbfam=\twelmsb
\textfont\msbfam=\tenmsb
\scriptfont\msbfam=\sevenmsb

\newtheorem{thm}{Theorem}[section]
\newtheorem{lemma}{Lemma}[section]
\newtheorem{defn}{Definition}[section]

\newcommand{\tr}{{\rm Tr}}
\newcommand{\ind}{{\rm Ind}}

\newcommand{\la}{\lambda}

\newcommand{\wi}{\widetilde}
\newcommand{\beq}{\begin{eqnarray}}
\newcommand{\eq}{\end{eqnarray}}
\newcommand{\bear}{\begin{array}}
\newcommand{\ar}{\end{array}}

\begin{document}

\renewcommand{\theequation}{\thesection.\arabic{equation}}
\setcounter{equation}{0}

\centerline{\Large{\bf On orbifold elliptic genus  }}
\vskip 7mm
\centerline{\bf Chongying DONG\footnote{Mathematics Department,
University of California,
Santa Cruz, CA 95064, U.S.A. (dong@math.ucsc.edu), Partially supported by NSF grant DMS-9987656 and a research grant from the Committee on Research, UC Santa Cruz.},
 Kefeng   LIU\footnote{Department of Mathematics, UCLA,CA 90095-1555,
USA (liu@math.ucla.edu). Partially supported by the Sloan
Fellowship and a NSF grant.}, Xiaonan
MA\footnote{Humboldt-Universit\"at zu Berlin, Institut f\"ur
Mathematik, Rudower Chaussee 25, 12489 Berlin, Germany.
(xiaonan@mathematik.hu-berlin.de). Partially supported by SFB
288.} }

\vskip 5mm {\bf 0 \hspace{.1in} Introduction}.
 \vspace{.1in}

  Elliptic genus was derived as
the partition function in quantum field theory \cite{W}.
Mathematically it is the beautiful combination of topology of
manifolds, index theory and modular forms (cf. \cite{L},
\cite{H}). Elliptic genus for smooth manifolds has been
well-studied. Recently, Borisov and Libgober ([3], [4]) proposed
some definitions of elliptic genus for certain singular spaces,
especially for complex orbifolds which is a global quotient $M/G$,
here the finite group $G$  acts holomorphically on complex variety
$M$. Similar definitions have been introduced by string theorists
in the 80s, in the study of orbifold string theory. One of their
guiding principles is the modular invariance. More recently
orbifold string theory has attracted the attention of geometers
and topologists. For example Chen and Ruan (cf. \cite{CR},\,
\cite{R})
  have defined orbifold
cohomology and orbifold quantum cohomology groups.

One of most important properties for elliptic genus is its
rigidity property under compact Lie group action. For smooth
manifolds, the rigidity and its generalizations have been well
studied. Since orbifold elliptic genus is the partition function
of orbifold string theory, it is natural to expect the rigidity
property for orbifold elliptic genus. Although the global
quotients from a very important class of orbifolds, many
interesting orbifolds are not global quotients. For example, most
of the Calabi-Yau hypersurfaces of weighted projective spaces are
not global quotients. In this paper we  define elliptic genus for
general orbifolds which generalizes  the definition of Borisov and
Libgober, and  prove their rigidity property. We actually
introduce the more general elliptic genus involving twisted
bundles and proved its rigidity. The idea of considering the
weights in the definition of orbifold elliptic genus comes from
\cite{BL} and  our proof of the K-theory version of Witten
rigidity theorems \cite[\S 4]{LMZ}. The proof of the rigidity is
again a combination of modular invariance and the index theory.
Only more complicated combinatorics are involved in the definition
and the proof.

This paper is organized as follows: In Sections 1 and 2 we review
the equivariant index theorem on orbifolds. We  define orbifold
elliptic genus and prove its rigidity for almost complex orbifolds
in Section 3.  Finally in Section 4 we introduce orbifold elliptic
genus for spin orbifolds, and will study its rigidity property on
a later occasion.

The authors would like to than Jian Zhou for many interesting
discussions regarding orbifold elliptic genus.

\section{ \normalsize Equivariant  index theorem for  spin orbifolds}
\setcounter{equation}{0}

In this section and the next section we recall the notations on
orbifolds, and explain the  equivariant index theorem for
orbifolds (cf. \cite[Chap. 14]{Du}, \cite{V}).

We first recall the definition of orbifolds, which are called
V-manifolds in \cite{K1}, \cite{Sa}.

We consider the pair $(G,V)$, here $V$ is a connected smooth
manifold, $G$ is a finite group acting smoothly and effectively on
$V$. A morphism $\Phi : (G,V)\to (G', V')$ is a family of open
embedding $\varphi: V\to V'$ satisfying:

{\bf i)} For  each $\varphi \in \Phi $, there is an injective group
 homomorphism
$\lambda_{\varphi}: G\rightarrow G' $ such that $\varphi$ is
$\lambda_{\varphi}$-equivariant.

{\bf ii)} For $g\in G', \varphi \in \Phi $, we define
$g\varphi :  V \rightarrow V'$
by $(g\varphi)(x) = g\varphi(x)$ for $x\in  V$.
If $(g\varphi)( V) \cap \varphi( V) \neq \phi$,
 then $g\in  \lambda_{\varphi}(G)$.

{\bf iii)} For $\varphi \in \Phi $, we have $\Phi = \{g\varphi, g\in G'\}$.\\

The morphism $\Phi$ induces a unique open embedding $i_{\Phi}: V/G
\to V'/G'$ of orbit spaces.

\begin{defn}  An orbifold $(X, {\cal U})$ is a paracompact Hausdorff space $X$ together with
a covering ${\cal U}$ of $X$ consisting of connected open subsets
such that

{\bf i)} For $U\in {\cal U}$, ${\cal V} (U)
=((G_U,\widetilde{U})\stackrel{\tau}{\rightarrow} U)$ is a
ramified covering $\widetilde{U}\rightarrow U$ giving an
identification $U \simeq \widetilde{U}/G_U$.

{\bf ii)} For $U, V\in {\cal U}, U\subset V$, there is a morphism
$\varphi_{VU}:(G_U,\widetilde{U})\rightarrow (G_V,\widetilde{V})$ that
 covers the
inclusion $U\subset V$.

{\bf iii)} For $U,V,W\in {\cal U}, U\subset V \subset W$, we have
$\varphi_{WU}=\varphi_{WV} \circ \varphi_{VU}$.
\end{defn}

In the above definition, we can replace  $(G,V)$ by a category of manifolds
 with an additional structure such as orientation, Riemannian metric
or complex structure.
We understand that the morphisms (and the groups) preserve the specified
 structure.
So we can define oriented, Riemannian or complex orbifolds.

\vspace{.1in}

 {\bf Remark}: Let $G$ be a compact Lie group and $M$ a smooth manifold
with a smooth $G$-action. We assume that the action of $G$ is
effective and  infinitesimally free. Then the quotient space $M/G$
is an orbifold. Reciprocally, any orbifold $X$ can be presented
this way. For example, let $O(X)$ be the total space of the
associated tangential orthonormal frame bundle. We know that
$O(X)$ is a smooth manifold and the action of the orthogonal group
$O(n)$ $(n= \dim X)$ is infinitesimally free on $O(X)$. The $X$ is
identified canonically with the orbifold $O(X)/O(n)$.\\

Let $X$ be an oriented orbifold, with singular set $\Sigma X.$ For
$x\in X$, there  exists a small neighbourhood $(G_x,
\widetilde{U}_x)\stackrel{\tau_x}{\to} U_x$ such that
$\widetilde{x}= \tau^{-1}_x(x) \in \widetilde{U}_x$
 is a fixed point of $G_x$. Such  $G_x$ is unique up to isomorphisms for each $x\in X$
 \cite[p468]{Sa}.
 Let (1), $(h_x^1), \cdots (h_x^{\rho_x})$ be all the conjugacy  classes
in $G_x$. Let $Z_{G_x}(h_x^j)$ be the centralizer of $h_x^j$ in
$G_x$. One also denotes  by $\widetilde{U}_x^{h_x^j}$ the  fixed
points of $h_x^j$ in $\widetilde{U}_x$. There is a natural
bijection
\begin{eqnarray}\label{0a1}
\{(y, (h_y^j)) |y\in U_x, j=1,\cdots \rho_y \} \\
\simeq \amalg_{j=1}^{\rho_x} \widetilde{U}_x^{h_x^j}/Z_{G_x}(h_x^j).
\nonumber
\end{eqnarray}
So we can define globally \cite[p77]{K1},
\begin{eqnarray}\label{0a2}
\widetilde{\Sigma X} =\{(x,(h_x^j))|x\in X, G_x \neq 1, j=1, \cdots \rho_x\}.
\end{eqnarray}
Then $\widetilde{\Sigma X }$ has a natural orbifold structure defined by
\begin{eqnarray}\label{0a3}
\Big \{(Z_{G_x}(h_x^j)/K_x^j, \widetilde{U}_x^{h_x^j})
\rightarrow  \widetilde{U}_x^{h_x^j}/Z_{G_x}(h_x^j) \Big \}_{(x, U_x,j)}.
\end{eqnarray}
Here $K_x^j$ is the kernel of the representation
$Z_{G_x}(h_x^j) \rightarrow$ Diffeo $(\widetilde{U}_x^{h_x^j})$.
The number $m=|K_x^j|$ is called the multiplicity of
$\widetilde{\Sigma X}$ in $X$ at $ (x, h_x^j)$.
Since the multiplicity is locally constant on $\widetilde{\Sigma X}$,
 we may assign the multiplicity  $m_i$ to each connected component
$X_i$ of $\widetilde{\Sigma X}$. In a sense $\widetilde{\Sigma X}$
is a resolution of singularities of $X$.

\begin{defn}
A mapping $\pi$ from an orbifold $X$ to an orbifold $X'$ is called
smooth if  for $x\in X, y=\pi(x)$, there exist  orbifold charts
$(G_x,\widetilde{U}_x), (G'_y,\widetilde{U'}_y) $ together with a
smooth mapping $\phi : \widetilde{U}_x \to \widetilde{U'}_y$ and a
homomorphism $\rho :  G_x \to G_y'$ such that $\phi$ is
$\rho$-equivariant and $\tau'_y\circ \phi  = \pi \circ \tau_x$.
\end{defn}

\begin{defn}
 An orbifold vector bundle $\xi$ over an orbifold
$(X,{\cal U})$ is defined as
follows: $\xi$ is an orbifold and for $U\in \cal U$,
$(G_U^{\xi}, \widetilde{p}_U: \widetilde{\xi}_U \rightarrow \widetilde{U})$ is
 a $G_U^{\xi}$-equivariant vector bundle and $(G_U^{\xi}, \widetilde{\xi}_U)$
(resp. $(G_U^{\xi}/K_U, \widetilde{U}), K_U= {\rm Ker}
(G_U^{\xi}\rightarrow {\rm Diffeo} (\widetilde{U})))$ is the
orbifold structure of $\xi$ (resp. $X$). In general, $G_U^{\xi}$
does  not  act effectively on $\widetilde{U}$, i.e. $K_U \neq \{
1\}$. If $G_U^\xi$ acts effectively on $\widetilde{U}$ for $U\in
\cal U$, we say  $\xi$ is a proper orbifold vector bundle.
\end{defn}

{\bf Remark}: Let $G$ be a compact Lie group acting effectivelty
and infinitesimally
freely on $M$. Then each $G$-equivariant bundle $E\to M$ defines a proper
orbifold vector bundle $E/G\to M/G$, and vice versa.\\

In the following, we will always denote by $(G_x,
 \widetilde{U}_x)$ $(x\in X)$ the orbifold chart as above.
For $h\in  G_x$, we have the following $h$-equivariant decomposition of
$T \widetilde{U}_x \otimes_{\sR} \bC$ as a real vector bundle on
$\widetilde{U}^h_x$,
\begin{eqnarray}
T \widetilde{U}_x  \otimes_\sR \bC = \bigoplus_{\lambda \in \bQ \cap
]0,1[} N_{\lambda (h)} \bigoplus T \widetilde{U}^h_x  \otimes_\sR \bC.
\end{eqnarray}
Here $N_{\lambda (h)}$ is the  complex vector bundle over
$\widetilde{U}^h_x$ with $h$ acting by $e^{2 \pi i \lambda }$ on
it. The complex conjugation provides a $\bC$ anti-linear
isomorphism between $N_{\lambda (h)}$ and
$\overline{N_{(1-\lambda) (h)}}$. If the order of  $h$ is even,
this produces a real structure on $N_{{1 \over 2}(h)}$, so this
bundle is the complexification of a real vector bundle
 $N_{{1 \over 2}(h)}^{\sR}$ on $\widetilde{U}^h_x$.
Thus, $T \widetilde{U}_x $ is isomorphic, as  a  real vector bundle, to
\begin{eqnarray} \label{1a1}
 T \widetilde{U}_x \simeq  \bigoplus_{\lambda \in \bQ \cap
]0,{1 \over 2}[} N_{\lambda (h)} \oplus N_{{1 \over 2}(h)}^{\sR}
\bigoplus T \widetilde{U}^h_x.
\end{eqnarray}
Note that $N_{\lambda (h)}$ (resp. $N_{{1 \over 2}(h)}^{\sR}$)
extends to complex (resp. real) vector bundle on
$\widetilde{\Sigma X}$. We will still denote them by $ N_{\lambda
(h)}$, $N_{{1 \over 2}(h)}^{\sR}$.

Assume that a compact Lie group $H$ acts differentially on $X$.
For $\gamma \in H$,
 let $X^\gamma = \{ x \in X, \gamma x = x \}$.
In the index theorem, we will use the following orbifold
as fixed point set of $\gamma$ which is a resolution
of singularities of $X^\gamma$ \cite[p180]{Du}.  For $x\in X^\gamma $,
 then on local chart $( G_x, \widetilde{U}_x)$,  $\gamma_{\widetilde{U}}$
acts on $\widetilde{U}_x$ as a linear map. The compatibility
condition for $\gamma_{\widetilde{U}}$ means that there exists an
automorphism $\alpha$ of $G_x$ such that for each $g\in G_x$,
$\gamma_{\widetilde{U}}\circ g \circ \gamma_{\widetilde{U}}^{-1} =
\alpha (g)$.  For $h\in G_x$, let $(h)_\gamma = \{ g h \alpha
(g)^{-1}; g\in G_x \}$ be the $\gamma$ conjugacy class in $G_x$.
Let
 \begin{eqnarray}\label{1a2}
\widehat{U}^{(h)_\gamma} _x = \{ (y,h_1)\in    \widetilde{U}_x\times G_x| (h_1 \circ
\gamma_   {\widetilde{U}}) (y)= y, h_1 \in (h)_\gamma\}.
\end{eqnarray}
Let $\widetilde{U}^{h\circ \gamma_{\widetilde{U}}}_x$ be the fixed point set
 of $h\circ \gamma_{\widetilde{U}}$ in $\widetilde{U}_x$, then
$\widetilde{U}^{h\circ \gamma_{\widetilde{U}}}_x$ is connected, and $x\in
\widetilde{U}^{h\circ \gamma_{\widetilde{U}}}_x$.

For $g\in G_x$, $g$ acts on $\widehat{U}^{(h)_\gamma} _x$ by the transformation
$$(y,h)\to (g(y), g\circ h\circ \alpha(g)^{-1}).$$
Indeed, if $(h\circ \gamma_{\widetilde{U}}) (y)=y$, as $\alpha(g)^{-1}
\circ  \gamma_{\widetilde{U}} =  \gamma_{\widetilde{U}}
\circ g^{-1} \circ\gamma_{\widetilde{U}}^{-1}  \circ\gamma_{\widetilde{U}}
=  \gamma_{\widetilde{U}}\circ  g^{-1}$, we know
 \begin{eqnarray}\label{0a5}
(gh\circ \alpha (g)^{-1} )\gamma_{\widetilde{U}}\circ g(y)
= gh\circ\gamma_{\widetilde{U}} (y)= g(y).
\end{eqnarray}

Let $Z^\gamma_{h,G_x} = \{ g\in G_x, gh\circ \gamma (g)^{-1} =h    \}$,
$K^\gamma_{h,G_x} = {\rm Ker}  \{ Z^\gamma_{h,G_x}  \to {\rm Diffeo}
 (\widetilde{U}^{h\circ \gamma_{\widetilde{U}} }_x) \}$. Then
 \begin{eqnarray}\label{0a6}
(Z^\gamma_{h,G_x}/K^\gamma_{h,G_x},
 \widetilde{U}^{h\circ \gamma_{\widetilde{U}}}_x)
\to \widetilde{U}^{h\circ \gamma_{\widetilde{U}}}_x/Z^\gamma_{h,G_x}
=\widehat{U}^{(h)_\gamma} _x/G_x\end{eqnarray}
defines an orbifold. We denote it by $\widetilde{X}^\gamma$.  Clearly,
$m(\widetilde{X}^\gamma) = |K^\gamma_{h,G_x} |$ is local constant on
$\widetilde{X}^\gamma$.

\begin{defn}\label{a5}
The oriented orbifold $X$ is spin if there exists 2-sheeted
covering of ${SO}(X)$ ($SO(X)$ is the oriented orthonormal frame
bundle of $TX$), such that for $U\in {\cal U}$, there exists a
principal Spin(n) bunlde ${\rm Spin}(\widetilde{U})$ on
$\widetilde{U}$, such that ${\rm Spin}(X)_{|U} \to {SO}(X)_{|U}$
is induced by ${\rm Spin}(\widetilde{U}) \to {SO}(\widetilde{U})$,
and ${\rm Spin}(\widetilde{U})$ also verifies the corresponding
compatible  condition.
\end{defn}

Then ${\rm spin}(X)$ is clearly a  smooth manifold.

Assume that orbifold  $X$ is spin. Let $h^{TX}$ be a metric on
$TX$ and $S(TX)= S^+(TX)\oplus S^-(TX)$ the corresponding orbifold
spinor bundle on $X$. Let $c(\cdot)$ be the Clifford action of
$TX$ on $S(TX)$. Let $\nabla^{S(TX)}$ be the connection on $S(TX)$
induced by the Levi-Civita connection $\nabla^{TX}$  on $TX$. Let
$W$ be a complex orbifold vector bundle on $X$. Let $\nabla^W$ be
a connection on $W$. Then $\nabla^{S(TX)\otimes W}= \
\nabla^{S(TX)}\otimes 1 + 1 \otimes \nabla^{W}$ is a connection on
$X$. Let $\Gamma (S^{\pm}(TX)\otimes W)$ be the set of $C^\infty$
sections of $S^{\pm}(TX)\otimes W$ on $X$. Let $D^X \otimes W$ be
the Dirac operator on $\Gamma (S^+(TX)\otimes W)$ to
 $\Gamma (S^-(TX)\otimes W)$, defined by
\begin{eqnarray}\label{0a7}
D^X \otimes W = c(e_i) \nabla_{e_i} ^{S^+(TX)\otimes W}.
\end{eqnarray}
Here $\{ e_i \}$ is an orthonrmal basis of $TX$.

Let $H$ be a compact Lie group.
 If $\gamma \in H$ acts on $X$ and lifts to ${\rm Spin}(X)$ and $W$.
Then $\nabla^{S(TX)}$ is $\gamma$ invariant and we can always find
$\gamma$ inavriant connection $\nabla^W$  on $W$. Note that
$D^X\otimes W$ is a $\gamma$ invariant elliptic operator on $X$.
For $x\in X$, let $K^W_x= {\rm Ker} (G_x^W \stackrel{\tau}{\to}
G_x)$. On $\widehat{U}^{(h)_ \gamma}$, let $N$ be the normal
bundle of
 $\widetilde{U}^{h\circ \gamma_{\widetilde{U}}}$ in $\widetilde{U}_x$. Let $W^0$ be the subbundle of
$W$ on $\widehat{U}_x$ which is $K^W_x$-invariant. Then $W^0$
extends to  a proper orbifold
vector bundle on $X$. We have the following decompositions:
\begin{eqnarray}\label{0a8}\begin{array}{l}
N= \oplus_{0< \theta < \pi} N_{\theta} \oplus N_{\pi},\\
W^0= \oplus_{0\leq  \theta < 2  \pi} W_\theta,
\end{array}\end{eqnarray}
where $N_{\theta}, W_\theta$ (resp. $N_\pi$) are complex (resp. real)
vector bundle on which $h\circ \gamma_{\widetilde{U}}$ acts as
multiplication by $e^{i \theta}$.  Then $\nabla^{TX}$ induces connection
$\nabla^{N_\theta}$ on $N_\theta$,
and $\nabla^{TX} = \oplus \nabla^{N_\theta}
\oplus \nabla^{TX^\gamma}$.
Let $R^W, R^{W^0},R^{N_\theta}, R^{TX^\gamma}$ be the curvatures of $\nabla^W, \nabla^{W^0},
 \nabla^{N_\theta},  \nabla^{TX^\gamma}$ ($\nabla^{W^0}$ is the connection on $W^0$ induced by $\nabla^W$).
\begin{defn} For $h\in G_x, g= h\circ \gamma_{\widetilde{U}}$, $ 0<\theta \leq \pi$, we write
\begin{eqnarray}\label{0a9}
\\
\begin{array}{l}
\ch_g(W, \nabla^W)= {1 \over |K_x^W|} \sum_{h_1\in  G_x^W, \tau(h_1)=h}
 \tr \Big [(h_1 \circ  \gamma_{\widetilde{U}}) \exp({ -R^{W}\over 2\pi i})\Big ]
=
\tr \Big [g \exp({ -R^{W^0}\over 2\pi i})\Big ],\\
\displaystyle{
\widehat{A}(T\widetilde{U}^g, \nabla^{T\widetilde{U}^g})
= {\det}^{1/2} \Big ({{i \over 4 \pi} R^{T\widetilde{U}^g}
\over \sinh({i \over 4 \pi} R^{T\widetilde{U}^g})}\Big )},\\
\displaystyle{
\widehat{A}_{\theta}(N_\theta, \nabla^{N_\theta})=
{ 1 \over     i^{{1 \over 2} \dim N_\theta} {\det}^{1 / 2} \Big (
1 - g \exp({i \over 2 \pi} R^{N_\theta}) \Big )},   }\\
\widehat{A}_g(N, \nabla^{N})= \Pi_{0<\theta \leq  \pi}
 \widehat{A}_{\theta}(N_\theta, \nabla^{N_\theta}).
\end{array} \nonumber
\end{eqnarray}
\end{defn}

If we denote by $\{x_j, -x_j \}$ $(j=1, \cdots, l)$  the Chern roots of
$N_\theta$, $T\widetilde{U}^g$  (Here we consider $N_\theta$ as a real vector bundle) such that $\Pi x_j$ defines the orientation of
$N_\theta$ and $T\widetilde{U}^g$, then
\begin{eqnarray}\label{0a10}\begin{array}{l}
\displaystyle{
\widehat{A}(T\widetilde{U}^g, \nabla^{T\widetilde{U}}) =\Pi_j { x_j \over 2}/ \sinh ({ x_j \over 2}),}\\
\displaystyle{
\widehat{A}_{\theta}(N_\theta, \nabla^{N_\theta})= 2^{-l} \Pi_{j=1}^l
{1 \over \sinh{1 \over 2} (x_j + i \theta)}
= \Pi_{j=1}^l { e^{{1 \over 2} (x_j + i \theta)}
\over e^{x_j + i \theta} -1}  .  }
\end{array}\end{eqnarray}

Recall that for $\gamma \in H$, the Lefschetz number ${\rm
Ind}_{\gamma} (D^X \otimes W)$, which is the
index of $D^X \otimes W$ if $\gamma =1,$  is defined by
\begin{eqnarray}
{\rm Ind}_{\gamma} (D^X \otimes W)= \tr \gamma _{|{\rm Ker} D^X\otimes W}-
\tr \gamma_{|{\rm Coker} D^X\otimes W}.
\end{eqnarray}
By using heat kernel, as in \cite[Th. 14.1]{Du}, we get
\begin{thm} \label{a1}  For $\gamma \in H$, we have the following equality:
\begin{eqnarray}\label{0a11}
{\rm Ind}_\gamma (D^X\otimes W)=\sum_{F\in \widetilde{X}^\gamma }
{1 \over m(F)} \int_{F} \alpha_F,
\end{eqnarray}
where
$\alpha_F$ is the
 characteristic class
\begin{eqnarray}\label{0a12}
\widehat{A}(T\widetilde{U}^{h\circ \gamma_{\widetilde{U}}},
\nabla^{T\widetilde{U}^{h\circ \gamma_{\widetilde{U}}}}) \Pi _{0< \theta \leq \pi}
 \widehat{A}_{\theta} (N_\theta, \nabla^{N_\theta})
{\rm ch}_{h\circ\gamma} (W, \nabla^{W})   \nonumber
\end{eqnarray}
on $\widetilde{U}_x^{h\circ \gamma_{\widetilde{U}}}.$
\end{thm}

Let $S^1$ act differentiably on $X$. Let $X^{S^1}= \{x \in X,
\gamma (x) =x, {\rm for \,  all } \, \gamma \in S^1  \}$. Let $V$
be the canonical basis of ${\rm Lie}(S^1) = \bR$. For $x\in X$,
let $V_X$ be the smooth vector field on $(G_x, \widetilde{U}_x)$
corresponding to $V$. Then $ V_X$ is
$G_x$-invariant \cite[p181]{Du}. We still denote by  $V_X$  the corresponding
smooth vector field on $X.$ We have $X^{S^1}= \{x \in X,V_X (x)=0\}$.

For $x\in X$, let $(1), \cdots (h^j_x), \cdots $ be the conjugacy
classes of $G_x$. Let $\widetilde{X}^{S^1} =  \{ (x, (h^j_x))| x\in
X^{S^1}, h^j_x \in G_x\}.$ Then $\widetilde{X}^{S^1}$ has a natural
orbifold structure defined by
\begin{eqnarray}\label{0a13}
\ \ \ \{Z_{G_x }(h^j_x)/K^{j,V }_x,
\widetilde{U}^{h^j_x}_V\!=\!\widetilde{U}^{h^j_x}_x \cap \{y\in
\widetilde{U}_x|V_X(y)\!=\!0 \} \} \to
(\widetilde{U}^{h^j_x}_V/Z_{G_x }(h^j_x),(h^j_x)),
\end{eqnarray}
where $K^{j,V }_x$ is the kernel of $Z_{G_x }(h^j_x) \to
{\rm Diffeo} \{\widetilde{U}^{h^j_x}_V\}$.

We have the following decomposition of smooth
vector bundles on $\widetilde{U}^{h}_V:$
\begin{eqnarray}\label{0a14}\begin{array}{l}
N_{\lambda (h)} = \oplus_j N_{\lambda, j},\\
N_{{1\over 2}(h)} ^\bR
= \oplus_{j>0} N_{{1\over 2},j}\oplus N^\bR_{{1\over 2},0},\\
T\widetilde{U}^{h} = \oplus_{j>0} N_{0,j} \oplus T\widetilde{U}^{h^j_x}_V,\\
W^0 = \oplus_{\la, j} W^0 _{\la, j}.
\end{array}\end{eqnarray}
Note that $N_{\lambda, j}, N_{{1\over 2},j}, N_{0,j} $  and $W^0 _{\la, j}$
extend to   complex vector bundles on $\widetilde{X}^{S^1}$, and
$\gamma =e^{2 \pi i t} \in S^1$
acts on them as multiplication by $e^{2 \pi i   j t}$. Also,
$N^\bR_{{1\over 2},0}$ and $ T\widetilde{U}^{h^j_x}_V$ extend
 to  real vector bundles on $\widetilde{X}^{S^1}$, and  $S^1$ acts on them as identity.  In fact, $T\wi{U}^h= T\widetilde{U}^{h^j_x}_V \oplus_{v\neq 0} N_{0,v,\sR}$, where $N_{0,v,\sR}$ denotes the undeling real bundle of the complex vector
 bundle $N_v$ on which $g\in S^1$ acts by multiplying by  $g^{v}$.
 Since we can choose either $N_v$
or $\overline{N}_v$ as the complex vector bundle for $N_{v, \sR}$, in
what follows, we always    assume $N_{{1\over 2},j}, N_{0,j}$ are zero if
$j<0$.

By (\ref{0a14}), for given $a\in \bC$, the eigenspace of
 $h\circ \gamma_{\widetilde{U}}$ with eigenvalue $a$ is equal to the
 sum of the above elements $N_{\lambda, j}$ such that
\begin{eqnarray}\label{0a15}
e^{2\pi i (\lambda + tj)} =a.
\end{eqnarray}
Let $A\subset \bR$ consist of $a\in \bR$ such that there exists
 $x \in X^{S^1}$, more than one  non-zero
 $N_{\lambda, j}$ on $\widetilde{U}^{h^j_x}_V$ are in the eigenspace
of $h\circ \gamma_{\widetilde{U}}$ with eigenvalue $e^{2 \pi i
a}$. As $X$ is compact, $A$ is a discrete set of $\bR$.

If $\gamma= e^{2 \pi i t}, t\in \bR\setminus  A$, then $\widetilde{X}^{\gamma}= \widetilde{X}^{S^1}$ by the construction. An immediate consequence
of Theorem \ref{a1} is the following.
\begin{thm} \label{a2} Under the condition of Theorem \ref{a1}, for $t\in \bR
\setminus A$, $\gamma = e^{2 \pi i t }$, we have
\begin{eqnarray}\label{0a16}
{\rm Ind}_\gamma (D^X\otimes W) = \sum_{F\in \widetilde{X}^{S^1} }
{1 \over m(F)} \int_{F} \alpha_F,
\end{eqnarray}
where $\alpha_F$ is the
 characteristic class
\begin{eqnarray}\label{0a17}
\widehat{A}(T\widetilde{U}^{h}_V, \nabla^{T\widetilde{U}^{h}_V} )
\sum_{\lambda, j} e^{2\pi i (\lambda + tj)} {\rm ch}
(W^0_{\lambda, j} , \nabla^{W^0})/\Pi_{\lambda, j}  i^{{1 \over 2}
\dim N_{\la,j}} {\det}^{1 / 2} \Big (
1 -  e^{2\pi i (\lambda + tj)}
\exp({i \over 2 \pi} R^{N_{\la,j}}) \Big )   \nonumber
\end{eqnarray}
on $\widetilde{U}^{h}_V.$
\end{thm}

\section{\normalsize Equivariant  index theorem for almost complex  orbifolds}
\setcounter{equation}{0}

If $X$ is an almost complex orbifold,  then
on the orbifold chart $(G_x, \widetilde{U}_x)$ for $x\in X,$
 we have the following
$h$-equivariant decomposition of  $T \widetilde{U}_x$ as  complex
vector bundles on  $\widetilde{U}^h_x$
\begin{eqnarray} \label{0a18}
 T \widetilde{U}_x \simeq  \bigoplus_{\lambda \in \bQ \cap
[0, 1[} N_{\lambda (h)}.
\end{eqnarray}
Here $N_{\lambda (h)}$ are complex vector bundles over
 $\widetilde{U}^h_x$ with $h$ acting by $e^{2 \pi i \lambda }$ on it,
and $N_{0(h)}$ is $T \widetilde{U}^h_x$. Again $N_{\lambda (h)}$
 extend to complex vector bundles on $\widetilde{\Sigma X}$.
We will still denote it by $ N_{\lambda (h)}$.

Let $F(x,h) = \sum _{\lambda }\lambda \dim_\sC N_{\lambda (h)}$ be the fermionic
shift, then $F: X\cup \widetilde{\Sigma X} \to \bQ$ is locally constant.
For  a connected component $X_i\subset X \cup \widetilde{\Sigma X}$ ,
we define $F(X_i)$ to be the values of $F$ on $X_i$.

Let $W$ be an orbifold complex vector bundle on $X$.
Let $D^X \otimes W$ be the Spin$^c$ Dirac operator on
$\Lambda (T^{*(0,1)} X)\otimes W$ \cite[Appendix D]{LaM}.

Let $H$ be a compact Lie group acting on $X$. We assume that the
action $H$ on  $X$  lifts on $W$, and  preserves the complex
structures of $TX$ and $W$. Now for  $\gamma\in H$, the
decomposition (\ref{0a8}) on $\widehat{U}^{(h)_ \gamma}_x$ also
preserves the complex structure of the normal bundle $N$. We
denote  by $R^{N}$ the curvature of
 $\nabla^{N}$ as  complex vector bundle. Then
\begin{eqnarray}\label{0a19}\begin{array}{l}
N= \oplus_{0< \theta <2 \pi} N_{\theta},\\
W^0= \oplus_{0\leq  \theta < 2  \pi} W_\theta.
\end{array}\end{eqnarray}
Here $N_{\theta}, W_\theta$  are complex
vector bundles on which $h\circ \gamma_{\widetilde{U}}$ acts as
multiplication by $e^{i \theta}$.
The following theorem is proved in  \cite[Th. 14.1]{Du}.
\begin{thm}\label{a3} Let
$${\rm Td}(T\widetilde{U}^{h\circ \gamma_{\widetilde{U}}},
\nabla^{T\widetilde{U}^{h\circ \gamma_{\widetilde{U}}}} )=\det
\left( {- R^{T\wi{U}^{h\circ \gamma_{\wi{U}}}} /2 i \pi \over 1-
\exp(- R^{T\wi{U}^{h\circ \gamma_{\wi{U}}}} /2 i \pi )}\right)$$
be the Chern-Weil Todd form of $ T\wi{U}^{h\circ
\gamma_{\wi{U}}}$. Then
                                       we have
\begin{eqnarray}\label{0a20}
{\rm Ind}_\gamma (D^X\otimes W) = \sum_{F\in \widetilde{X}^\gamma }
{1 \over m(F)} \int_{F} \alpha_F.
\end{eqnarray}
Here on $\widetilde{U}^{h\circ \gamma_{\widetilde{U}}}$,
$\alpha_F$ is the
 characteristic class
\begin{eqnarray}\label{0a21}
{\rm Td}(T\widetilde{U}^{h\circ \gamma_{\widetilde{U}}},
\nabla^{T\widetilde{U}^{h\circ \gamma_{\widetilde{U}}}} ) {\rm
ch}_{h\circ \gamma}  (W, \nabla^{W})/ \det (1- (h\circ \gamma)
\exp ({i \over 2 \pi} R^{N} )).
 \nonumber
\end{eqnarray}
\end{thm}

If $H=S^1$,   on $\widetilde{U}^{h}_V$ as in (\ref{0a13}), we have the
following decomposition of complex vector bundles,
\begin{eqnarray}\label{0a22}\begin{array}{l}
N_{\lambda (h)} = \oplus_j N_{\lambda (h), j},\\
T\widetilde{U}^{h} = \oplus N_{0,j} \oplus T\widetilde{U}^{h^j_x}_V.
\end{array}\end{eqnarray}
Here $N_{\lambda (h), j}, N_{0,j} $ extend to   complex vector
bundles on $\widetilde{X}^{S^1}$, and $\gamma =e^{2 \pi i t}\in
S^1$ acts on them as multiplication by $e^{2 \pi i   j t}$.

By  Theorem \ref{a3}, we get,
\begin{thm} \label{a4} Under the condition of Theorem \ref{a3}, for $t\in \bR
\setminus A$, $\gamma = e^{2 \pi i t }$, we have
\begin{eqnarray}\label{0a21}
{\rm Ind}_\gamma (D^X\otimes W) = \sum_{F\in \widetilde{X}^{S^1} }
{1 \over m(F)} \int_{F} \alpha_F.
\end{eqnarray}
Here on $\widetilde{U}^{h}_V$, $\alpha_F$ is the
 characteristic class
\begin{eqnarray}\label{0a22}
{\rm Td}(T\widetilde{U}^{h}_V, \nabla^{T\widetilde{U}^{h}_V} )
\sum_{\lambda, j} e^{2\pi i (\lambda + tj)} {\rm ch}
(W^0_{\lambda, j} , \nabla^{W^0}) /\Pi_{\lambda, j}   {\det} \Big
( 1 -  e^{2\pi i (\lambda + tj)}\exp({i \over 2 \pi}
R^{N_{\lambda, j}}) \Big ).
           \nonumber
\end{eqnarray}
\end{thm}

Note that we can also get Theorems \ref{a2} and  \ref{a4} from
\cite[Theorem 1]{V}.

\section{ \normalsize Elliptic genus for almost complex orbifolds}

\setcounter{equation}{0}

In this section, we define the elliptic genus for a general
 almost complex  orbifold  and prove its rigidity property.
We are using  the setting of Section 2.

For $\tau \in \bH = \{ \tau \in \bC; {\rm Im} \tau >0\}$,
 $q= e^{ 2\pi i \tau}$, $t\in \bC$,  let
\begin{eqnarray} \label{0b1}\begin{array}{l}
\theta (t, \tau)=c(q)q^{1/8} 2 \sin (\pi t )
\Pi_{k=1}^\infty (1 - q^{k} e^{2 \pi i t})
\Pi_{k=1}^\infty (1 - q^{k} e ^{-2 \pi i t}).
\end{array}\end{eqnarray}
be the classical Jacobi theta function \cite{Ch},
where $c(q)= \Pi_{k=1}^\infty (1 - q^{k} )$.   Set
\begin{eqnarray}
\theta'(0,\tau) = {\partial \theta (\cdot, \tau) \over \partial t}|_{t=0}.
\end{eqnarray}
Recall  the following transformation formulas for the
theta-functions \cite{Ch}:
\begin{eqnarray}\label{0b2}\begin{array}{l}
\theta (t+1, \tau) = -\theta (t,\tau), \qquad
\theta ( t+ \tau, \tau)= - q^{-1/2} e^{- 2 \pi i t} \theta (t, \tau),\\
\theta ({t \over \tau}, - {1 \over \tau})= {1 \over i} \sqrt{\tau \over i}
 e^{\pi i t^2 \over \tau} \theta (t,\tau), \quad
 \theta (t, \tau+1) = e^{ \pi i \over 4} \theta (t, \tau).
\end{array}\end{eqnarray}

For a complex or  real vector bundle $F$ on a manifold $X$, let
\begin{eqnarray}\label{0b3}\begin{array}{l}
\Sym_q (F) = 1 + q F + q^2 \Sym^2 F + \cdots,\\
\Lambda_q (F) = 1 + qF + q^2 \Lambda^2 F + \cdots,
\end{array}\end{eqnarray}
be the symmetric and the exterior power  operations $F,$
respectively.

Let $X$ be an almost complex orbifold, and $\dim _\sC X = l$.
In this Section, all vector bundles are complex vector bundles.
Let $W$ be a proper orbifold complex  vector bundle on $X$,
and $\dim_\sC W = m$.
Then $W^0$ in (\ref{0a19}) is $W$. Now for the vector bundle $W$, the
fermionic shift $F(X_i,W) =\sum _\la \la \dim W_\la $ is well define
on each connected component $X_i\subset X \cup \wi{\Sigma X}$.
For $x\in X$, $y=e^{2 \pi i z}$,
we use the orbifold  chart  $( G_x, \widetilde{U}_x)$.
For $h\in G_x$, by (\ref{0a18}), we define  on $\widetilde{U}_x^h$,
\begin{multline} \label{0b4}
\Theta^z_{q, x,  (h)} (TX) =  \bigotimes_{\lambda \in \bQ \cap [0,1[ }
\Big (  \bigotimes_{k=1}^{\infty}
\Big (\Lambda_{-y^{-1}q^{k-1 + \lambda (h)}} N^*_{\lambda (h) }
\otimes
 \Lambda_{-yq^{k - \lambda (h)}} N_{\lambda (h) } \Big ) \Big ) \\
\bigotimes_{\lambda \in \bQ \cap ]0,1[ } \Big ( \bigotimes_{k=1}^{\infty}
\Big (\Sym_{q^{k-1 + \lambda (h)}} N^*_{\lambda (h) }
\otimes
 \Sym_{q^{k -\lambda (h)}} N_{\lambda (h) } \Big )\Big ) \\
\bigotimes_{k=1}^{\infty}\Big (\Sym_{q^{k}} N^*_{0 }
\otimes\Sym_{q^{k }} N_{0 } \Big ),
\end{multline}
\begin{multline}
\Theta_{q, x,  (h)}^z (TX|W) =  \bigotimes_{\lambda \in \bQ \cap [0,1[ }
\Big (  \bigotimes_{k=1}^{\infty}
\Big (\Lambda_{-y^{-1}q^{k-1 + \lambda (h)}} W^*_{\lambda (h) }
\otimes
 \Lambda_{-yq^{k - \lambda (h)}} W_{\lambda (h) } \Big ) \Big ) \\
\bigotimes_{\lambda \in \bQ \cap ]0,1[ } \Big ( \bigotimes_{k=1}^{\infty}
\Big (\Sym_{q^{k-1 + \lambda (h)}} N^*_{\lambda (h) }
\otimes
 \Sym_{q^{k -\lambda (h)}} N_{\lambda (h) } \Big )\Big ) \\
\bigotimes_{k=1}^{\infty}\Big (\Sym_{q^{k}} N^*_{0 } \otimes
 \Sym_{q^{k }} N_{0 } \Big ). \nonumber
\end{multline}
It is easy to  verify that
each coefficient of $q^a$ $(a\in \bQ)$
in $\Theta_{q, x,  (h)}^z (TX)$, $\Theta_{q, x,  (h)}^z (TX|W)$
 defines an orbifold vector bundle on $\widetilde{\Sigma X}$. We will denote
 it by $\Theta_{q, X_i  }^z (TX)$, $\Theta_{q, X_i}^z (TX|W)$ on the connected component $X_i$
of $X \cup \widetilde{\Sigma X}$. It is the usual Witten element
on $X$ (see \cite{H} and \cite{W}),
\begin{eqnarray}\label{0b5} \qquad
\Theta_q^z (TX) =  \bigotimes_{k=1}^{\infty} \Big
(\Lambda_{-y^{-1}q^{k-1}} T^* X \otimes \Lambda_{-yq^{k }} TX\Big
) \bigotimes_{k=1}^{\infty} \Big (\Sym_{q^{k}} T^* X \otimes
\Sym_{q^k} TX  \Big ).
\end{eqnarray}

\begin{defn} \label{b1}
The orbifold elliptic genus of $X$ is defined to be
\begin{eqnarray}\label{0b6}
F(y,q) = y^{ l \over 2} \sum_{X_i \subset X \cup \widetilde{\Sigma X}}
y^{-F(X_i)} \ind (D^{X_i}\otimes  \Theta_{q, X_i}^z (TX)).
\end{eqnarray}
More generally, we define the orbifold elliptic genus associated
to $W$ as
\begin{eqnarray}\label{0b61}
F(y,q, W) = y^{ m \over 2} \sum_{X_i \subset X \cup \widetilde{\Sigma X}}
y^{-F(X_i,W)} \ind (D^{X_i}\otimes  \Theta_{q, X_i}^z (TX|W)).
\end{eqnarray}
\end{defn}

If $X$ is a global quotient $M/G$ where the action
of finite group $G$ on almost complex manifold $M$ preserves its
complex structure, the equation (\ref{0b6}) coincides with
\cite[Definition 4.1]{BL}.

We  next prove that the orbifold elliptic genus is rigid for $S^1$
action on $X$.

Let $S^1$ act on $X$, preserving  the complex structure on $TX$, and lifting on $W$.
Naturally, we define the Lefschetz number for $\gamma \in S^1$,
\begin{eqnarray}\label{0b7}\qquad
F_\gamma (y,q,W) = y^{m \over 2} \sum_{X_i \subset X \cup \widetilde{\Sigma X}}
y^{-F(X_i,W)} \ind_\gamma (D^{X_i}\otimes  \Theta_{q, X_i}^z (TX|W))
\end{eqnarray}

Let $P$ be a compact manifold acted infinitesimally freely by a compact Lie
group $G, $ and $X=P/G$ the corresponding orbifold.
We still denote $W$ the corresponding vector bundle on $P$ for $W$.
 Then $K_X= \det (T^{(1,0)} X)$, $K_W=\det W$ are  naturally induced by   complex line bundles
on $P$, we still denote it by $K_X$, $K_W$. We may also consider $K_X$, $K_W$
as an orbifold line bundle on $X$.

Recall that the equivariant cohomology group $H^*_{S^1} (P, \bZ)$
of $P$ is defined by
\begin{eqnarray}
H^*_{S^1} (P, \bZ)= H^*(P \times_{S^1} ES^1, \bZ).
\end{eqnarray}
where $ES^1$ is the universal $S^1$-principal bundle over
the  classifying space $BS^1$ of $S^1$.
So $H^*_{S^1} (P, \bZ)$ is a module over $H^*(BS^1, \bZ)$ induced by the
projection $\overline{\pi} : P\times _{S^1} ES^1\to BS^1$.
Let $p_1(W)_{S^1}, p_1(TX)_{S^1} \in H^*_{S^1} (P, \bZ)$ be the equivariant
 first Pontrjagin classes of $W$ and $TX$ respectively. Also recall that
\begin{eqnarray}
H^*(BS^1, \bZ)= \bZ [[u]]
\end{eqnarray}
with $u$ a generator of degree $2$.

Recall from \cite[Theorem B]{Liu2} that for smooth manifold $X$
one needs the conditions \beq \label{1b12} p_1(W-TX)_{S^1} = 0, \
\  c_1(W-TX)_{S^1} =0. \eq for the rigidity theorem.

Note that if the connected component $X_i$ of $X\cup
\widetilde{\Sigma X}$ is  defined by $(\wi{U}_x^h, Z_{G_x}(h))$,
for $g\in Z_{G_x}(h)$, set 
$U^{h,g}_V=\{ b\in \widetilde{U}_x| hb=g b= b, V_X(b)=0\}$.
Then the connected component $X_{ik}$ of $\wi{X^{S^1}_i}$  is
defined by $({U}^{h,g}_V, Z_{Z_{G_x}(h) }(g))$. We have the
following decomposition of complex vector bundles on $ U^{h,g}_V$,
\begin{align}\label{1b13}
T\widetilde{U}_x = \sum_{\lambda(h), \lambda(g) \in \sQ \cap
[0,1[, v \in \sZ} N_{\lambda(h), \lambda(g), v },\\ \nonumber W =
\sum_{\lambda(h), \lambda(g) \in \sQ \cap [0,1[, v \in \sZ}
W_{\lambda(h), \lambda(g), v }.\nonumber
\end{align}
where $g,h \in G_x$ (resp. $\gamma = e^{2 \pi i t}\in S^1$) act on
$N_{\lambda(h), \lambda(g), v }, W_{\lambda(h), \lambda(g), v }$
as multiplication by $e^{2 \pi i \lambda(g)}$, $e^{2 \pi i
\lambda(h)}$ (resp. $e^{2 \pi i vt})$. Let $2 \pi i
x^j_{\lambda(h), \lambda(g), v }$, $2 \pi i w^j_{\lambda(h),
\lambda(g), v }$ be the formal Chern roots of $N_{\lambda(h),
\lambda(g), v }$, $W_{\lambda(h), \lambda(g), v }$ respectively.
To simplify the notation, we will omit the superscript $j$.

Then $N_{\la (h), \la(g), v}$, $W_{\lambda(h), \lambda(g), v }$ extend to
 orbifold vector bundles on $X_{ik}$.  Now the natural generalization of
(\ref{1b12})  for orbifold is the following: there exists $n\in
\bN$, such that on  each
 connected component $X_{ik}$  of $\wi{X^{S^1}_i}$,
\begin{align}\label{1b14}
&\sum_{\la (h), \la(g), v,j} \Big [ (w_{\la (h), \la(g), v}^j +   \la (g) -\tau \la (h)  + vu)^2  \\ \nonumber
&\hspace*{10mm}-  (x_{\la (h), \la(g), v}^j +   \la (g) - \tau \la (h)  + vu)^2    \Big ]
 = n\overline{\pi}^* u^2 \in H^*(X_{ik}, \bQ) [\tau,t],
\nonumber
\end{align}
\begin{align}\label{1b15}
& \sum_{\la (h), \la(g), v,j} (w_{\la (h), \la(g), v}^j +   \la (g) - \tau \la (h)+ vu)\\
&\hspace*{10mm}-    \sum_{\la (h), \la(g), v,j} (x_{\la (h), \la(g), v}^j +   \la (g) - \tau \la (h)+ v u ) =0 \in H^*(X_{ik}, \bQ) [\tau,t].\nonumber
\end{align}

\begin{thm}\label{b2} Assume that
 $S^1$ acts on $P$ which induces the $S^1$-action
on $X$, and lifts to $W$, and  $ c_1 (W)=0 \mod\!N$ in $H^*(P,
\bZ)$ for some $1<N\in \bN.$  Also assume that  there exists $n\in
\bN$ such that equations (\ref{1b14}) and  (\ref{1b15}) hold. Then
 for any  $N$-th root of unity $ y = e^{2 \pi i z}$ we have

i) If $n=0$, then
 $F_\gamma (y,q, W)$
is constant on $\gamma \in S^1$.

ii) If $n<0$, then
 $F_\gamma (y,q, W)=0$.
\end{thm}

Note that, in case $W=TX$, the condition (\ref{1b14}), (\ref{1b15})
 are automatic, and
as a consequence we get the rigidity and vanishing theorems for
the usual orbifold elliptic genus $F(y, q)$. In particular we know
that for a Calabi-Yau almost complex manifold $X$, $F(y, q)$ is
rigid for any $y$.

\vspace{.2in}

 $Proof$: Using Theorem \ref{a4}, for $\gamma = e^{2
\pi i t}, t\in \bR \setminus A$, $ y=e^{2\pi i z}, q=e^{2\pi i
\tau}$, we get
\begin{eqnarray}\label{0b8}
F_\gamma (y,q,W) = y^{m \over 2} \sum_{X_i \subset X \cup \widetilde{\Sigma X}}
y^{-F(X_i,W)} \sum_{F \subset  \widetilde{X}_i^{S^1}}
 {1 \over m(F)} \int_{F} \alpha_F,
 \end{eqnarray}
Recall that $V_X$ is the smooth vector field generated
by $S^1$-action on $X$. For $x\in X$, take the orbifold chart
  $( G_x, \widetilde{U}_x)$. 
If $X_i \subset X\cup \widetilde{\Sigma X}$ is  represented by
$\widetilde{U}_x^h/Z_{G_x}(h)$
on $\widetilde{U}_x$ as in (\ref{0a3}), the normal bundle
$N_{X_i,g,V}= N_{\widetilde{U}_x^h/U^{h,g}_V}$ of $U^{h,g}_V$ in
$\widetilde{U}_x^h$ extends to an orbifold
vector bundle on $\widetilde{X}_i^{S^1}$.
By Theorem \ref{a4}, the contribution of the chart $(G_x, \widetilde{U}_x)$
for $\ind_\gamma (D^{X_i}\otimes \Theta_{q, X_i}^z  (TX|W))$ is
\begin{eqnarray}\label{0b9}
{1 \over |Z_{G_x}(h)|}\sum_{gh=hg, g\in G_x}\int_{U^{h,g}_V}
{{\rm Td}(TU^{h,g}_V) {\rm ch}_{g\circ \gamma}{(\Theta_{q, x,  (h)}^z (TX|W))}
  \over \det (1- (g\circ \gamma) e ^{-{1 \over 2 \pi i} R^{N_{X_i, g,V}}})}.
\end{eqnarray}
 Let
\begin{align}\label{1b18}
N_v = \sum_{\lambda(h), \lambda(g) \in \sQ \cap [0,1[}
N_{\lambda(h), \lambda(g), v },\\ \nonumber
W_v = \sum_{\lambda(h), \lambda(g) \in \sQ \cap [0,1[}
W_{\lambda(h), \lambda(g), v }.\nonumber
\end{align}
As the vector field $V_X$ commutes with the action of  $G_x$,
 $N_v$ and $ W_v$ extend to vector bundles on $ \widetilde{X}_i^{S^1}$.
  So the contribution of the chart $( G_x, \widetilde{U}_x)$
for $F_\gamma (y,q,W)$ is
\begin{eqnarray}\label{0b10}
\\
\begin{array}{l}
\displaystyle{{1 \over |G_x|}\sum_{gh=hg;  g,h \in G_x}y^{{m \over 2}-F(X_i,W)}
\int_{U^{h,g}_V}  {{\rm Td}(TU^{h,g}_V)
{\rm ch}_{g\circ \gamma}{\Theta_{q, x,  (h)}^z  (TX|W)}
\over \det (1- (g\circ \gamma) e ^
{-{1 \over 2 \pi i} R^{N_{X_i,g}}}) }  }\\
\displaystyle{={1 \over |G_x|}\sum_{gh=hg;  g,h \in G_x}\int_{U^{h,g}_V}
( 2 \pi i x_{0(h),0(g),0}) } y^{{m \over 2}-F(X_i,W)} \\
\displaystyle{\hspace*{10mm}\Pi_{\lambda(g), v } {1\over 1-e^{2\pi i (x_{0(h), \lambda(g), v }+ \lambda(g)+t v)}}
   \Pi_{\lambda(h)>0, \lambda(g), v } } \\
\displaystyle{\hspace*{10mm}
 \Big (\Pi_{k=1}^{\infty}  {(1-y^{-1}q^{k-1 + \lambda (h)}
e^{2\pi i (-w_{\lambda(h), \lambda(g), v }- \lambda(g)- t v)})
(1-yq^{k - \lambda (h)}
e^{2\pi i (w_{\lambda(h), \lambda(g), v }+ \lambda(g)+ t v)})
\over (1-q^{k-1 + \lambda (h)}
e^{2\pi i (-x_{\lambda(h), \lambda(g), v }-\lambda(g)- t v)})
(1-q^{k - \lambda (h)}
e^{2\pi i (x_{\lambda(h), \lambda(g), v }+\lambda(g)+ t v)})}  \Big )}    \\
\displaystyle{\hspace*{10mm}
\Pi_{\lambda(g), v } \Pi_{k=1}^{\infty}  {(1-y^{-1}q^{k-1}
e^{2\pi i (-w_{0(h), \lambda(g), v }- \lambda(g)- t v)})
(1-yq^{k}
e^{2\pi i (w_{0(h), \lambda(g), v }+ \lambda(g)+ t v)})
\over (1-q^{k}
e^{2\pi i (-x_{0(h), \lambda(g), v }-\lambda(g)- t v)})
(1-q^{k}
e^{2\pi i (x_{0(h), \lambda(g), v }+\lambda(g)+ t v)})} } \\
\displaystyle{
= (i^{-1} c(q) q^{1/8} )^{l-m} {1 \over |G_x|}\sum_{gh=hg;  g,h \in G_x}\int_{U^{h,g}_V} ( 2 \pi i  x_{0(h),0(g),0}) }\\
\displaystyle{\hspace*{10mm} {\Pi_{\lambda(h), \lambda(g), v,j }
(\theta (w^j_{\lambda(h), \lambda(g), v }
+ \lambda(g)-\tau \lambda(h) + z + tv, \tau) e^{-2 \pi i z \lambda (h)})\over
\Pi_{\lambda(h), \lambda(g), v,j }\theta (x^j_{\lambda(h), \lambda(g), v }
+ \lambda(g)-\tau \lambda(h) + tv, \tau)} }.
\end{array} \nonumber
\end{eqnarray}
To get the last equality of (\ref{0b10}), we use  (\ref{0b1}), (\ref{1b15}).

We will consider $F_\gamma (y,q,W)$ as a function of $(t,z,\tau)$,
we can extend it to a meromorphic function on $\bC \times \bH
\times \bC$. From now on, we denote $(i^{-1} c(q) q^{1/8} )^{m-l}
{\theta'(0,\tau)^l \over  \theta (z,\tau)^m} F_\gamma (y,q,W)$  by
$F(t,\tau,z)$. We also denote  by $F_\gamma
(t,\tau,W)_{|\widetilde{U}_x}$ the function defined by
(\ref{0b10}).
 We set
 \beq F(t,\tau,z)_{|\widetilde{U}_x} =(i^{-1}
c(q) q^{1/8} )^{m-l} {\theta'(0,\tau)^l \over  \theta (z,\tau)^m}
F_\gamma (t,\tau,W)_{|\widetilde{U}_x}. \eq

Now, the equation (\ref{1b14}) implies the equalities
\begin{align}\label{1b21}
& \sum w_{\la (h), \la (g),v}^2 -\sum x _{\la (h), \la (g),v}^2 =0,
\ \  \sum v \   w_{\la (h), \la (g),v} - \sum v \   x_{\la (h), \la (g),v} =0,
\\ \nonumber
&\sum_{\la (h), \la (g)}  \la (h) \la (g) (\dim W_{\la (h), \la (g)}
 - \dim N_{\la (h), \la (g)}) =0, \\ \nonumber
&\sum_{\la (h), \la (g),v} \la (g) v \  (\dim W_{\la (h), \la (g), v}
- \dim N_{\la (h), \la (g), v}) =0, \\ \nonumber
&\sum   \la (g)^2  (\dim W_{\la (g)}- \dim N_{\la (g)})=0,
\ \  \sum    v^2 (\dim W_v -\dim N_v)= n.  \nonumber
\end{align}

By (\ref{0b1}),  for $a,b \in 2 \bZ $, $ k\in \bN$,
\begin{eqnarray}\label{1b22}
\theta (x + k(t + a \tau + b), \tau) = e^{- \pi i (2k a x + 2 k^2 a t
+ k^2 a^2 \tau)} \theta (x + kt , \tau).
\end{eqnarray}

As $c_1(K_X) =0 \mod\!N$ in $H^* (P,\bZ)$, by the same argument as
\cite[\S 8]{H} or \cite[Lemma 2.1, Remark 2.6]{LMZ},
 $\sum_u v \dim N_v \mod\!N$ is constant on each connected
component of $X$.

By (\ref{0b10}), (\ref{1b21}), (\ref{1b22}), we know for $a,b
\in 2  \bZ$,
\begin{eqnarray}\label{0b14}
F(t+a \tau +b,\tau,z)
=  e^{-2 \pi i z a \sum_v  v \dim N_v } F(t,\tau,z).
\end{eqnarray}
For $A= \left ( \begin{array}{l} a \quad b\\
c \quad d \end{array} \right ) \in SL_2 (\bZ)$,  we define its modular
transformation on $\bC \times \bH$ by
\begin{eqnarray}\begin{array}{l}
\displaystyle{
A(t, \tau) = \left ( { t \over c \tau + d}, {a \tau + b \over c \tau + d}
\right ). }
\end{array}\end{eqnarray}
By (\ref{0b10}),
 under the action $A= \left ( \begin{array}{l} a \quad b\\
c \quad d \end{array} \right ) \in SL_2 (\bZ)$,
we have
\begin{align}\label{1b24}
&{\theta(A(t_1,\tau))
\over \theta(A(t_2,\tau))} =
e^{\pi i c  (t^2_1-t^2_2) / c \tau + d }
{\theta(t_1, \tau) \over \theta(t_2, \tau)},  \\ \nonumber
&    {\theta' (A(0,\tau))
\over \theta(A(t,\tau))} = (c \tau +d)
e^{- \pi i c  t^2 / c \tau + d }
{\theta ' (0, \tau) \over \theta(t, \tau)}.    \nonumber
\end{align}

 For $g,h\in G_x$, by looking at the degree $2 \dim_\sC U^{h,g}_V$ part, that is the $\dim_\sC U^{h,g}_V$-th
homogeneous terms of the polynomials in $x, w$'s,
 on both sides of the following equation, we get
\begin{eqnarray}\label{0b17} \qquad \begin{array}{l}
 \displaystyle{  \int_{U_V^{h,g}}
( x_{0(h),0(g),0})\Pi_{\lambda(h), \lambda(g), v }
e^{2 \pi i c z  w_{\lambda(h), \lambda(g), v }
(c \tau + d) }  }\\
 \displaystyle{
 {\theta (w_{\lambda(h), \lambda(g), v }(c \tau +d)
+ \lambda(g)(c \tau +d)-(a \tau + b) \lambda(h) + z(c \tau +d) + tv,\tau)\over
\theta (x_{\lambda(h), \lambda(g), v }(c \tau +d)
+ (c \tau +d)\lambda(g)-(a \tau + b)\lambda(h) + tv,\tau)} }\\
 \displaystyle{
= \int_{U_V^{h,g}} ( x_{0(h),0(g),0})\Pi_{\lambda(h), \lambda(g), v }
e^{2 \pi i c z w_{\lambda(h), \lambda(g), v }} }  \\
 \displaystyle{  {\theta (w_{\lambda(h), \lambda(g), v }
+ \lambda(g)(c \tau +d)-(a \tau + b) \lambda(h) + z(c \tau +d) + tv,\tau)\over
\theta (x_{\lambda(h), \lambda(g), v }
+ (c \tau +d)\lambda(g)-(a \tau + b)\lambda(h) + tv,\tau)} }.
\end{array}\end{eqnarray}
By (\ref{0b2}), (\ref{0b10}), (\ref{1b21}),  (\ref{1b24}) and
(\ref{0b17}), we easily derive the following identity:
\begin{multline} \label{0b18}
 F( A(t,\tau),z)_{|\widetilde{U}_x}
 =
{1 \over |G_x|} (c \tau + d)^l e^ { \pi i c n t^2 /(c\tau + d)}
{\theta'(0,\tau)^l \over   \theta (z(c\tau + d),\tau)^m}\\
\sum_{gh=hg;  g,h \in G_x}\int_{U^{g,h}_V}
( 2 \pi i x_{0(h),0(g),0})  \\
 \times    \Big \{ \Pi_{\lambda(h), \lambda(g), v }
(e^{2 \pi i c z  (w_{\lambda(h), \lambda(g), v }
  +(\lambda (g)-{a \tau + b \over c \tau + d}\lambda (h))
(c \tau + d) + t v  ) } e^{-2 \pi i z \la (h)})\Big \}  \\
 \times  {\Pi_{\lambda(h), \lambda(g), v }  \theta (w_{\lambda(h), \lambda(g), v }
+ \lambda(g)(c \tau +d)-(a \tau + b) \lambda(h) + z(c \tau +d) + tv,\tau)
 \over
\Pi_{\lambda(h), \lambda(g), v }  \theta (x_{\lambda(h), \lambda(g), v }
+ (c \tau +d)\lambda(g)-(a \tau + b)\lambda(h) + tv,\tau) }.
\end{multline}
 By (\ref{0b2}), (\ref{1b15}), (\ref{1b21}), (\ref{1b22}), (\ref{0b18}),
 we have
\begin{align} \label{1b28}
&(c \tau + d)^{-l} e^ {-\pi i c n t^2 /(c\tau + d)}
{\theta (z(c\tau + d),\tau)^m\over \theta'(0,\tau)^l}
F( A(t,\tau),z)_{|\widetilde{U}_x} \\ \nonumber
&\hspace*{3mm}= {1 \over |G_x|}\sum_{gh=hg;  g,h \in G_x}\int_{U^{h,g}_V}
 (2 \pi i  x_{0(h),0(g),0}) \\ \nonumber
& \hspace*{10mm}\Pi_{\lambda(h), \lambda(g), v,j } \Big \{ e^{2 \pi i c z\Big (d\lambda (g)-b\lambda(h) + \tau (c\lambda(g)-a \lambda(h))\Big ) }  e^{2 \pi i c z (w^j_{\lambda(h), \lambda(g), v }+tv)} \\ \nonumber
& \hspace*{15mm}\times e^{-2 \pi i z (c\lambda(g)-a \lambda(h) + \lambda(g^{-c}h^{a})) (c\tau +d)}e^{-2 \pi i z \lambda (h)}) \Big  \}
\\   \nonumber
& \hspace*{10mm}
 \times {\Pi_{\lambda(h), \lambda(g), v }\theta (w_{\lambda(h), \lambda(g), v }
+ \lambda(g^dh^{-b})-\tau  \lambda(g^{-c}h^{a}) + z(c \tau +d) + tv,\tau)\over
\Pi_{\lambda(h), \lambda(g), v } \theta (x_{\lambda(h), \lambda(g), v }
+ \lambda(g^dh^{-b})-\tau  \lambda(g^{-c}h^{a}) + tv,\tau)} \\ \nonumber
&\hspace*{3mm}= {1  \over |G_x|}\sum_{gh=hg;  g,h \in G_x}
\int_{U^{h,g}_V} ( 2 \pi i x_{0(h),0(g),0})
 \Pi_{\lambda(h), \lambda(g), v }   \Big (e^{2 \pi i c z (w_{\lambda(h),
  \lambda(g), v } +tv)} e^{-2\pi i z (c \tau +d)\lambda(g^{-c}h^{a})} \Big )
\\ \nonumber
 &\hspace*{10mm}
 {\Pi_{\lambda(h), \lambda(g), v } \theta (w_{\lambda(h), \lambda(g), v}
+ \lambda(g^dh^{-b})-\tau  \lambda(g^{-c}h^{a}) + z(c \tau +d) + tv,\tau)
 \over
\Pi_{\lambda(h), \lambda(g), v } \theta (x_{\lambda(h), \lambda(g), v }
+ \lambda(g^dh^{-b})-\tau  \lambda(g^{-c}h^{a}) + tv,\tau)  }. \nonumber
\end{align}
 Recall that $c_1(K_W)= 0 \mod \!N$ where $K_W={\mbox{det}}\, W$, this implies that the line bundle $K^{cz}_W$
 is well defined on $P$, also as an orbifold vector bundle on $X$. Let
\begin{multline}\label{1b29}
F^A (t,\tau,z) = (i^{-1} c(q) q^{1/8} )^{m-l}{\theta'(0,\tau)^l \over   \theta (z(c\tau + d),\tau)^m}\\
\sum_{X_i \subset X \cup \wi{X}}  y^{{m\over 2} - F(X_i,W)}{\rm Ind}_\gamma \Big (D \otimes K^{cz}_W \otimes
\Theta^{(c\tau +d)z}_{q,X_i} (TX|W)\Big ).
\end{multline}
Now, observe that as  $A= \left ( \begin{array}{l} a \quad b\\
c \quad d \end{array} \right ) \in SL_2 (\bZ)$,
when $g,h$ run through all pair of $\sum_{gh=hg, g,h\in G_x}$,
then $g^{-c}h^{a}, g^d h^{-b}$ run through all pair of
$\sum_{gh=hg, g,h\in G_x}$. Then by (\ref{1b28}), (\ref{1b29})
\begin{eqnarray}
F(A(t,\tau),z)=(c \tau + d)^l e^ { \pi i c n t^2 /(c\tau + d)}
F^A (t,\tau,z).
\end{eqnarray}

The following lemma implies that the index theory comes in to cancel
part of the poles of the functions $F$.
\begin{lemma}\label{b3}
The function $F^A (t,\tau,z)$ is holomorphic in $(t, \tau)$
for $(t, \tau)\in \bR \times \bH$.
\end{lemma}

The proof of the above Lemma is the same as the proof of
\cite[Lemma 1.3]{Liu2} or
 \cite[Lemma 2.3]{LMa1}.

Now, we  return to the proof of Theorem \ref{b2}. Note that the possible polar divisors of $F$ in $\bC
\times \bH$ are of the form $t= {k \over j} (c\tau + d)$ with $k,
c, d, j$ integers and $(c,d)=1$ or $c=1$ and $d=0$.

We can always find integers $a, b$ such that $ad-bc = 1$, and
consider the matrix
$A=\left ( \begin{array}{l}d\quad -b\\
-c \quad a
\end{array} \right ) \in SL_2(\bZ)$.
\begin{eqnarray}\qquad \begin{array}{l}
\displaystyle{
 F^A(t, \tau,z) =(-c \tau + a)^{-l} e^ {\pi i c n t^2 /(-c\tau + a)}
 F \Big (A(t,\tau), z \Big )}
\end{array}\end{eqnarray}

Now, if $t= {k \over j} (c\tau + d)$ is a polar divisor of $F (t,\tau,z)$,
then one polar divisor of $F^A (t,\tau, z)$ is given by
\begin{eqnarray}
{t \over -c \tau + a}= {k \over j}
 \Big ( c {d \tau -b \over -c \tau + a} + d \Big ),
\end{eqnarray}
which exactly gives $t = k/j$. This contradicts Lemma \ref{b3}, and completes the proof of Theorem \ref{b2}. \hfill $\blacksquare$

\section{\normalsize Elliptic genus
for spin orbifolds }

\setcounter{equation}{0}

We are following the setting of  Section 1.

For $\tau \in \bH = \{ \tau \in \bC; {\rm Im} \tau >0\}$,
 $q= e^{ 2\pi i \tau}$, let
\begin{eqnarray}\label{0c1}\begin{array}{l}
\theta_3(v, \tau)=c(q)\Pi_{k=1}^\infty (1 + q^{k-1/2} e^{2 \pi i v})
\Pi_{k=1}^\infty (1 + q^{k-1/2} e ^{-2 \pi i v}).
\end{array}\end{eqnarray}
be the other three classical Jacobi theta-functions \cite{Ch},
where $c(q)= \Pi_{k=1}^\infty (1 - q^{k} )$.

Let $X$ be a compact orbifold, $\dim_\sR X = 2n$. We assume that $X$ and $\widetilde{\Sigma X}$
are spin in the sense of Definition \ref {a5}.
 For $x\in X$,
taking the orbifold  chart  $(G_x, \widetilde{U}_x)$.
By (\ref{1a1}), for $h\in G_x$, we define  on $\widetilde{U}_x^h$
\begin{eqnarray}\label{0c2}\begin{array}{l}
\Theta'_{q, x,  (h)} (TX) =
 \otimes_{\lambda \in \bQ \cap ]0,{1\over 2}[ }
\Big ( \otimes_{k=1}^{\infty}
\Big (\Lambda_{q^{k-1 + \lambda (h)}} N^*_{\lambda (h) }
\otimes
 \Lambda_{q^{k - \lambda (h)}} N_{\lambda (h) } \Big ) \\
\hspace*{25mm}\otimes_{k=1}^{\infty}
\Big (\Sym_{q^{k-1 + \lambda (h)}} N^*_{\lambda (h) } \otimes
 \Sym_{q^{k - \lambda (h)}} N_{\lambda (h) } \Big ) \Big )\\
\hspace*{10mm}\otimes_{k=1}^{\infty}
\Big (\Lambda_{q^{k-{1\over 2}}} N_{{1\over 2}(h)}^\bR  \otimes
\Sym_{q^{k-{1\over 2}}} N_{{1\over 2}(h)}^\bR \Big )
\otimes \otimes_{k=1}^{\infty}  (\Lambda_{q^k} (T\widetilde{U}_x^h) \otimes
\Sym_{q^k} (T\widetilde{U}_x^h)).
\end{array}\end{eqnarray}

It is easy to verify that each coefficient of $q^a$ $(a\in \bQ)$
in $\Theta'_{q, x,  (h)} (TX)$
 defines an orbifold vector bundle on $X\cup \widetilde{\Sigma X}$.
We  denote  as $\Theta'_{q, X_i  } (TX)$ on the connected component
 $X_i$ of $ X\cup \widetilde{\Sigma X}$. Especially,
 $\Theta'_{q, X} (TX)$
is the usual Witten elements on $X$
\begin{eqnarray}\label{0c3}\begin{array}{l}
\Theta'_q (TX) =  \bigotimes_{k=1}^{\infty}
\Big ( \Lambda_{q^{k }}( TX) \otimes\Sym_{q^k} (TX)  \Big ).
\end{array}\end{eqnarray}

We propose the following definition for the elliptic genus on spin
orbifold:

\begin{defn}
The orbifold elliptic genus of $X$ is
\begin{eqnarray}\label{0c4}
F(q) = \sum_{X_i \subset X \cup \widetilde{\Sigma X}}
{\rm Ind} \Big (D^{X_i}\otimes (S^+(TX_i)\oplus S^-(TX_i))  \otimes \Theta'_{q, X_i} (TX)\Big ).
\end{eqnarray}
\end{defn}

Let $S^1$ act on $X$ and preserve the spin structure of $X \cup
\widetilde{\Sigma X}$. Naturally, we define the Lefschetz number for
 $\gamma\in S^1$
\begin{eqnarray}\label{0c5}\qquad
F_{d_s,\gamma} (q) = \sum_{X_i \subset X \cup \widetilde{\Sigma X}}
 {\rm Ind}_\gamma\Big  (D^{X_i}\otimes (S^+(TX_i)\oplus S^-(TX_i)) \otimes  \Theta_{q, X_i} (TX)\Big )
\end{eqnarray}

On local chart $(G_x, \widetilde{U}_x)$, for $h,g\in G_x$, $gh=hg$,
 then as in (\ref{1a1}), on $\widetilde{U}_x^h$, we have
\begin{eqnarray}\label{1c10}\begin{array}{l}
T\widetilde{U}_x = N_0 \oplus_{\lambda(h)\in ]0,{1 \over 2}[} N_{\lambda(h)}
\oplus N^\sR_{{1\over 2}(h)}.
\end{array}\end{eqnarray}
here $h$ acts on the real vector bundles $N_0=T\widetilde{U}_x^h$,
$N_{{1\over 2}(h)}^\sR$ as multiplication by $1, e^{\pi i}$. $h$ acts
on complex vector bundles $N_{\lambda(h)}$ as multiplication by
$e^{2 \pi i \lambda(h)}$. Now on $\widetilde{U}_x^{h,g}$, the
fixed point set of $g$ on  $\widetilde{U}_x^{h}$, we have the
following decomposition
\begin{eqnarray}\label{0c10}\begin{array}{l}
N_{\lambda(h)} = \oplus_{\lambda(g)\in [0,{1 \over 2}]} N_{\lambda(h), \lambda(g)}
\oplus_{\lambda(g)\in ]0,{1 \over 2}[} N_{\lambda(h), 1-\lambda(g)}
\quad {\rm for} \quad 0< \lambda(h)< {1\over 2},\\
N_0= \oplus_{\lambda(g)\in [0,{1 \over 2}]} N_{0,\lambda(g)},\\
N^\sR_{{1\over 2}(h)}= \oplus_{\lambda(g)\in [0,{1 \over 2}]}
N_{{1 \over 2}(h), \lambda(g)}.
\end{array}\end{eqnarray}
Here $N_{\lambda(h), \lambda(g)}$ $(\lambda(h), \lambda(g)\in
\{0, {1\over 2}\})$ are real vector bundles on $\widetilde{U}_x^{h,g}$.
And  $N_{\lambda(h), \lambda(g)}$ $(\lambda(h)$ or $ \lambda(g)$ not in $\{0, {1\over 2}\})$
 are complex  vector bundles on $\widetilde{U}_x^{h,g}$.
$h,g$ act on $N_{\lambda(h), \lambda(g)}$ as multiplication by
$e^{2 \pi i \lambda(h)}$, $e^{2 \pi i \lambda(g)}$ respectively.
Again $N_{\lambda(h), \lambda(g)}$ extends to a vector bundle on
 $\widetilde{X}_i^{S^1}$.

For $f(x)$ a holomorphic function, we denote by
$f(x_{\la(h),\la(g)})(N_{\la(h),\la(g)}) = $ $ \Pi_j f(x_{\la(h),\la(g),j})$,
the symmetric polynomial which gives characteristic class of
$N_{\la(h),\la(g)}$.
 we get that  the contribution of the chart $(
G_x, \widetilde{U}_x)$ for $F_{d_s, \gamma} (q)$  is
\begin{eqnarray} \label{0c15}\begin{array}{l}
\displaystyle{F_{d_s} (t,\tau)_{|\widetilde{U}_x} =
{i^{-n} \over |G_x|}\sum_{gh=hg; g,h\in G_x}\int_{U^{h,g}_V}
\Big ( 2 \pi i  x_{0(h),0(g)}  {\theta_1 (x_{0(h),0(g)},\tau) \over
\theta (x_{0(h),0(g)},\tau)}\Big ) (N_{0(h),0(g)} )}\\
\displaystyle{\hspace*{10mm}\Pi_{\stackrel{\lambda(h)\in \{0,{1\over 2}\}, 0\leq \lambda(g)\leq {1\over 2}}{(\lambda(h), \lambda(g))\neq (0,0)} }
{\theta_1 (x_{\lambda(h), \lambda(g)}
+ \lambda(g)-\tau \lambda(h) ,\tau)\over
\theta (x_{\lambda(h), \lambda(g)}
+ \lambda(g)-\tau \lambda(h),\tau)}(N_{\lambda(h), \lambda(g)})}\\
\displaystyle{\hspace*{15mm}\Pi_{0<\lambda(h)< {1\over 2}} \Big [ \Pi_{ 0\leq
\lambda(g)\leq {1\over 2}}{\theta_1 (x_{\lambda(h), \lambda(g) }
+ \lambda(g)-\tau \lambda(h),\tau)\over
\theta (x_{\lambda(h), \lambda(g)}
+ \lambda(g)-\tau \lambda(h),\tau)}(N_{\lambda(h), \lambda(g)})}\\
\displaystyle{\hspace*{18mm} \Pi_{0<\lambda(g)< {1\over 2}}  {\theta_1 (x_{\lambda(h), 1-\lambda(g)}
- \lambda(g)-\tau \lambda(h),\tau)\over
\theta (x_{\lambda(h),1- \lambda(g)}
- \lambda(g)-\tau \lambda(h),\tau)}(N_{\lambda(h),1- \lambda(g)})\Big ]}
\end{array}\end{eqnarray}

We plan to return to the study of their rigidity and vanishing
properties on a later occasion.

\end{document}